\documentclass{article}
\usepackage{graphicx} %
\usepackage{tcolorbox}
\usepackage[numbers]{natbib}

\title{Mathematical Folklore}
\author{Govind Menon\thanks{Division of Applied Mathematics, Brown University, 182 George St., Providence, RI 02912.} \thanks{Supported by the NSF grant DMS 2407055 and the Erik Ellentuck Fellowship at the Institute for Advanced Study, Princeton. } and Akshay Venkatesh\thanks{Institute for Advanced Study, 1 Einstein Drive, Princeton NJ.}}
\date{\today}

\begin{document}

\maketitle

What is human in mathematics? This is the fundamental question at the heart of an experimental seminar, {\em Mathematical Folklore}, that we organized in 2025/2026 
at the Institute for Advanced Study (IAS) in Princeton. We were inspired, in part, by the Frankfurt Mathematical History seminar, led by Max Dehn—and beautifully memorialized by Siegel \cite{siegel_frankfurt_seminar}. Yet we wanted to not only study great works of mathematics but also experience its folk aspects—the myriad tendrils by which it draws nourishment from surrounding culture.   And so we drew no line between scientific and general writings, reading dense mathematical prose side by side with friendly correspondence, essays on culture, and theological discussions.

  In organizing the seminar, we sought an intellectual basis to think about the turmoil that currently surrounds our discipline: If mathematics can be mechanized, what should we relegate to machine, and what do we keep---what {\em can} we keep---for humans? What {\em is} the basic value of our subject, and how do we articulate it in a world of fractured attention, when the academy and the structures that support it seem to be in disarray? But we wanted, too, a space for types of reflection and community that are often squeezed out by our professional commitments. We wanted to suggest a slower, more deliberate, type of scholarship whereby we read, think and discuss questions that are outside of mathematics proper: Where does this esoteric topic truly come from? Why do we prove things at all? How does society shape mathematics? How has mathematics shaped society? We often think of our discipline as a collection of formal statements—but how do human beings relate to that abstract world?

 We felt that there was no better way to approach these issues than by examining the historical development and evolution of the subject. All of us, implicitly, are in constant dialogue with our past; when we teach a class or write a book, we reorganize a historical process of knowledge creation into a linearly ordered chain of definitions and theorems. Such abstraction allows us to share knowledge with ease, but it obscures the view of mathematics as a living organism in constant evolution. It excludes a great part of those traditions and cultures that give mathematics its vitality,
 and it is to these that the title of our seminar refers.

We sought, in our seminar, to reverse our usual pedagogical practice---to recontextualize, to informalize---to let mathematicians speak for themselves, both in the past and in the present. Listening to voices from the seventeenth and nineteenth centuries, we quickly found that the roots of the subject are far messier---far more {\em human}---than our usual narratives suggest: they lie as much in astrology as astronomy, and as much in philosophy as physics. By forgetting these origins, we forget, too, part of ourselves; we forget that our subject was shaped by accidents and mistakes, choices and values, as much as by the needs of science and the laws of logic. And in this collective loss of memory, we lose awareness of the possibility of making different choices for the future, limiting our dreams.

This article presents reflections on our experience: how the seminar was organized, lessons we learned from it, and some seedlings of a vision for humanism and mathematics. And while what we have written represents our view as organizers, it reflects and integrates many perspectives that we learned from participants in the seminar over the course of the year; we thank all of them for a thrilling journey. We would also like to thank an anonymous referee and Siobhan Roberts
for helpful suggestions and edits.

\section{Organizing the seminar}
We will share a few details of how our seminar was organized—who participated, the structure of each session, the curriculum and how it was chosen, and how it might be adapted to other contexts.

\subsection{Who came?}

We started the year with trepidation, not knowing if anyone would participate. But we were consistently surprised by the broad enthusiasm in our community for the project. In the end, our group was diverse, including mathematicians at different career stages, graduate students and postdocs with varied specializations, historians of science, science journalists and writers, and occasional visitors to the IAS. We cannot emphasize enough how much the seminar was enriched by   participants outside the community of professional mathematicians—both through the perspectives they brought to the discussions and through the various readings they curated.

\subsection{The structure of the seminar}
 Most mathematics seminars consist of a speaker explaining ideas to an audience. We wanted to create a seminar that included more voices, but also ensure that the discussions remained rigorous and rooted in the readings, with sufficient pauses to think. We sat around a table, but with a blackboard close at hand—discussions frequently included a brief presentation or spontaneous explanation in appreciation of an ingenious argument or technique, such as Klein’s surprising use of projective geometry as the foundation for non-Euclidean geometry, or various propositions on conics in Newton’s {\em Principia}. Although not all the audience would engage with the more technical material, there was no shortage of it for those who were interested—topics ranged from statistical mechanics to aperiodicity to hyperbolic trigonometry.
 
We also needed to engage an audience with widely different backgrounds. To this end, we split our reading into (frequently difficult) core material and easier supplementary material such as correspondence from archives, commentaries, historical analysis and popular writing.

Each discussion had four designated readers—the two of us and two volunteers—who would often meet beforehand to talk through ideas. The seminar began with 20-30 minutes of commentary by the readers, followed by an hour of open discussion that would often continue into lunchtime.

Not everyone read everything, nor was everyone present at each discussion. However, we found that the group quickly developed an atmosphere of trust and warmth; in a typical session, most participants would speak.

\subsection{The curriculum}
The readings for our seminar are described in the Appendix. The first term traces a path through 19th-century mathematics leading to the crisis of foundations and the creation of the computer. In the second term, armed with a better sense of how mathematicians of the 19th century looked at the past, we decided to follow their gaze back to the beginnings of the Scientific Revolution.

While we set our goals in broad outline—How did the Erlangen program arise? Why was the molecular theory of matter controversial?—and chose some seminal articles as anchors, the topics evolved in response to the discussions. We did not plan more than three discussions ahead, nor did we work through the readings completely before assigning them. This choice reflects our experience as academics, but also the fact that we were designing the seminar on the fly. It had the advantage of bringing to life aspects of the past that we (mistakenly) thought we understood. For example, reading Klein inspired us to return to Lobachevsky, and reading Newton inspired us to return to Kepler; and in both cases, reading the earlier work changed the way we viewed the later.

In general, as the year went on and the rapport amongst the group grew, we took a more relaxed approach as organizers. We allowed ourselves many sessions that were not tightly connected to the larger arc—sessions covering amateur mathematics, the “laws of thought,” the role of AI in current research, and  women in mathematics. These themes were suggested by seminar participants, who curated the readings with us.

In general, the task of reading itself was harder than we expected and (lacking experience and driven by curiosity) we assigned too much. A modern mathematics text uses a standardized expository style with definitions, theorems, and proofs. It is primarily technical, with a few guarded remarks of a contextual nature. The papers we read were not of this form. They were uneven in format and length, written with unfamiliar language and notation, and expressed in a tone that was often conceptual and philosophical. Thus, considerable effort was needed to undo our familiar reading habits.

Many classical sources, such as the work of Kepler and Newton, have many translations and are easily accessible in the public domain. In less explored territory, we found chatbots (mindfully used)  an invaluable tool for translation and contextual research. And, in a time when many of our younger colleagues may never stand in a mathematical library, we discovered the joy of deep dives into archival material, such as a trove of correspondence between Norman Steenrod and his undergraduate advisor Raymond Wilder.

\subsection{History and heritage}

The approaches of historians and mathematicians to the history of mathematics have often been at odds (see, for instance, Grattan-Guinness’s distinction \cite{GG} between “history” and “heritage,” or a more detailed taxonomy in \cite{Fried}). But it is difficult to locate our seminar along this spectrum, for our goals were personal; we sought inspiration and rejuvenation as much as knowledge. We were inspired by certain histories of mathematics, and followed a somewhat conventional historical arc, but we focussed on original texts rather than secondary sources. And while we inevitably approached those texts as practising mathematicians, we evaluated ideas in the rough, placed restrictions on the methods we were allowed to use, and separated the historical and logical development of mathematics when we could.

\subsection{Adapting the seminar}
 
 While several aspects of this experimental undertaking reflect the particular freedom provided by the Institute for Advanced Study, the seminar’s basic structure is very flexible and accessible---indeed, our seminar emerged from the liberal arts curriculum, in the form of a seminar at Brown University for undergraduates designed by Lindsay Caplan and the first author. Others, too, are organizing seminars that integrate mathematics and metamathematics—for example, Aravind Asok’s Folklore Seminar at USC \cite{asok_seminar}, and an online seminar Math(s), Philosophy, History run by the Free Computing Lab. And in a more technical vein, Chris Skinner and the second author have organized multiple seminars on the history of number theory that have proved quite popular with graduate students.
 
These experiences suggest that the structure of a seminar of this type can be adapted to the interests of a wide range of participants. Mathematics is blessed with a deep textual tradition, which can be interrogated in many ways, at a cadence suited to the audience. Our curriculum should not be seen as a list of readings to be replicated, but just one re-imagining of the past among many.

One element of our project, however, strikes us as fixed. It is necessary to grapple with original sources, and it is necessary to slow down. To pass too quickly over missteps and flaws, to dismiss the unfamiliar tone of past work, is to miss the first flicker of a new imagination. Progress in mathematics is contingent on many circumstances, but none is more important than the free will of the scholars who came before us.

\section{Reflections} \label{Reflections}
We will describe a few themes that kept recurring over the course of the year, along with some salient quotations from our readings. They each illustrate, in their own ways, that mathematics is, and has always been, rooted in very human concerns.

\subsection{The formal and the intuitive}
   \begin{quote}
{\em The viewpoint of the formalist must lead to the conviction ... [that] Not to the mathematician, but to the psychologist, belongs the task of explaining why we believe in certain systems of symbolic logic and not in others, in particular why we are averse to the so-called contradictory systems...} ---Brouwer \cite{brouwer_intuitionism}
  \end{quote}

 Brouwer’s essay, an important salvo in the crisis of foundations, puts sharp focus on a collection of ideas that we would encounter again and again. To paraphrase Hermann Weyl, it seems that the angel and devil of intuition and formalism (we leave the reader to pick which is which) have battled for the soul of mathematics in every generation.
 
For Leibniz, as today, a formal language enabled the definitive resolution of disputes. But it frequently has served a purpose quite different to Brouwer’s caricature above: as a tool for studying our own minds. Hilbert’s earliest formalist writing,  {\em Grundlagen der Geometrie},  begins by quoting Kant, and suggests that the axiomatic analysis is “tantamount to the logical analysis of our intuition of space.” Half a century earlier, Boole had sought similarly to understand “the fundamental laws of reasoning in the symbolic language of a Calculus.” Later, the emerging category theory of Eilenberg and Maclane would “provide a technical background for the intuitive notion of naturality.”

And on the other side, as our quote above suggests—intuitionism for Brouwer reflected not just a set of axioms but a view that formalism separated mathematics from the human in a profound way. This theme could be seen in many other readings, from Newton’s deep antipathy to Descartes’ algebraic calculus, to the question, as we entered the  automated theorem proving era in the 1950s, of how computer reasoning would relate to human intuition. Dick’s historical study \cite{dick_aftermath} of Argonne’s AURA program of the 1970s and 1980s provides us with a harbinger of what we may expect, at much larger scale, in the next decade:

\begin{quote}
{\em ... Overbeek acknowledged that the kind of intuition needed to collaborate with AURA successfully was not simply logical or mathematical intuition ...    it was technological and computational intuition.}
\end{quote}

Seeing these arguments  endure over centuries  we wonder: is the crisis of foundations truly settled?

\subsection{Mathematics as a social system}
\begin{quote}

{\em Do I think things are better? Of course. There were few women in the Greek world except goddesses, and less to be said of the age of enlightenment. Here in the US, women got the vote in 1920, and Title IX was passed some fifty years later in 1972. We had the second wave of feminism in the 1960’s, the founding of AWM in 1971, the evolution of the AMS from an openly misogynous organization to one proclaiming diversity, and more than two or three token women professors in most good mathematics departments.}   ---Uhlenbeck \cite{uhlenbeck_glass_ceiling}
\end{quote}

 Mathematics today has emerged from waves of collective effort, driven by social and political events that lie beyond our control. The Cold War and the collapse of the Soviet Union, the creation of funding agencies and national labs, marches for civil rights, decolonization in the Global South, new growth in Asia---all have changed the fabric of  academic mathematics.  In the other direction, the majority of our PhD students go on to careers outside the academy, bringing an advanced mathematical training to a great variety of professions. How do we begin to understand the complexity of this story?
 
Our seminar tugged at a few of these strands, even without conscious effort— for the act of reimagining a text also produces an imagination of the underlying social context. It was particularly vivid to see how these forces shaped the lives of our participants. In one session, Karen Uhlenbeck reflected on the fragility and unevenness of the social mobility that had brought her and a pioneering cohort of women into mathematics. In another, we were struck by the depth of enthusiasm amongst our peers for the tradition of mathematics sometimes called “recreational,” but that we would like to call {\em folk}---a form of mathematics, passed on through puzzles and columns and games, that still draws many of us into the discipline. As Doris Schattschneider writes in her celebration of
Marjorie Rice and Martin Gardner:

\begin{quote}
{\em  The mind and spirit are the forte of all such amateurs --- the intense spirit of inquiry and the keen perception of all they encounter. No formal education provides these gifts. Mere lack of a mathematical degree separates these ``amateurs'' from the ``professionals''. Yet their dauntless curiosity and ingenious methods make them true mathematicians.} ---Schattschneider  \cite{Schattschneider}
\end{quote}

 The relationship between mathematics and social power is complicated, and several discussions in our community at this juncture have been deeply divisive. The nature of this community and the boundaries of mathematics are themselves in flux; we need to reflect on our autonomy with greater rigor. A folklore seminar provides, at least, a space for an examination of fundamental values which, in its way, is the start of an antidote to tribalism.

\subsection{The many languages of mathematics}
 \begin{quote}
     {\em Just as languages like Greek or Sanskrit are historical facts and not absolute logical necessities,
it is only reasonable to assume that logics and mathematics are similarly historical, accidental forms of expression.} ---von Neumann \cite{vonneumann_computer_brain}
 \end{quote}
 As our year proceeded, we were struck indeed by the sheer variety of languages in which we have spoken of mathematics.
 
Thus, for example, the analysis of planetary motion—what would now be treated as a nonlinear differential equation and analyzed by means of calculus or symplectic dynamics—is overcome, by Newton, through a {\em tour de force} of Euclidean geometry. Or the ``monster of multilinear algebra'' that is the Riemann curvature tensor is described in a single paragraph of natural language by its inventor, making use of just one equation. Both texts are, without doubt, difficult to read—but they have a richness, texture and accessibility that seem to have been lost as we adopted a more abstracted language. At times we felt, with regret, that we have traded sourdough for Wonder Bread.

Even the 20th-century texts we studied, such as Wiener’s analysis of Brownian motion or the early mathematical writing of Eilenberg and Steenrod, diverge significantly from modern presentations. Returning to history allows us to perceive the slow but inexorable process of ongoing translation in mathematics. And it makes one wonder: Who will read our papers in a hundred years and what will they make of them?

\subsection{Mathematics and the senses}
\begin{quote}
{\em How does instinctive knowledge originate and what are its contents? Everything we analyze in nature imprints itself {\em uncomprehended\/} and {\em unanalyzed\/} in our percepts and ideas, which, then in their turn, mimic the processes of nature in their most striking features. In these accumulated experiences we possess a treasure-store which is ever close at hand of which only the smallest portion is embodied in articulate thought\/} ---Mach ~\cite{mach_mechanics}
\end{quote}

As we discussed Riemann and Helmholtz, and observed the process of formalization that followed, we were struck by the care with which the 19th-century scientists approached the distinction between our inner world and the nature of experimental data. Helmholtz’s work on the axioms of geometry was strongly influenced by his work on physiology. Both Riemann and Hilbert engaged carefully with Kant, critiquing his views on {\em a priori} knowledge of space. Mach’s empiricist philosophy of science, grounded in his work on the psychology of perception, influenced Einstein’s early work on relativity, especially in the idea of measurement by observers. When viewed in tandem, our readings allowed us to see the revolutions in geometry and physics as being guided by philosophies of science and experiments on perception.

Vision, in particular, was a source of mathematical thought throughout our readings. Descartes, Leibniz and Newton all studied optics (and Newton experimented rather fearlessly on his own eyes). Helmholtz was an expert on vision, and Riemann’s paper \cite{riemann_hypothesen} offers two examples of notions corresponding to continuous manifolds: colour, and the positions of perceived objects. We also read Mumford---much later and driven in part by the study of computer vision---proposing that the foundations of mathematics be based on probability rather than logic.

\subsection{Faith}

\begin{quote}
    {\em Accordingly the movements of the heavens are nothing except a certain
everlasting polyphony (intelligible, not audible) with dissonant tunings...
it is no longer a surprise that man, the ape of his Creator, should
finally have discovered the art of singing polyphonically.\/} ---Kepler~\cite{kepler_harmonice}
\end{quote}

We found persistent links between mathematics and religious faith in many of our readings. From our modern standpoint, it is tempting to disregard Newton’s analysis of biblical prophecy or Kepler’s quest for celestial polyphony---but these mystical ideas are manifested, both subtly and concretely, in their mathematical work. Kepler, for example, spent decades seeking rational ratios in the numerology of planetary motion, guided by his ``irrational'' beliefs in heavenly harmony. Working through the Principia allowed us to appreciate Newton’s mastery of technique; but reading his work on the occult allowed us to enter what Einstein called his {\em geistige Werksttat}. A comparative analysis of his writings reveals a careful fidelity to textual analysis and a fervent desire to recover {\em prisca sapientia}, the lost wisdom of the ancients.  In his closeted heresy and zealotry, in his deep mistrust of the power of institutions, we see a complex human being, driven by extreme ambition to understand all there is to know about the world. Mathematics and terror have never felt as close.

While less overt, religious influences were still perceptible in several readings from the 19th century. The theme of ancient wisdom is explicit in  Maxwell’s work on molecules~\cite{maxwell_molecules}, and it was tempting to read  Riemann’s emphasis on simplicity, clarity and immediacy of the mathematical experience as a reflection of his Lutheran faith.

 Coming to the 20th century, we read the
    correspondence of the famous algebraist Andr{\'e} Weil with
    his even more famous sister, the  mystic and philosopher Simone Weil,
     who challenged  Andr{\'e} to account for the meaning of mathematics. 
     Constrasting the modern notion of mathematics-as-game with 
     her image of the ancient Greek tradition, she wrote: 

\begin{quote}
    {\em Purity of soul was their only concern; “imitating God” was its secret; the study of mathematics helped to imitate God insofar as one saw the universe as subject to mathematical laws, which made the geometer an imitator of the supreme legislator ---}S. Weil \cite{weil_correspondence}
\end{quote}
We went back to some of the ancient  sources, where we found similar
sentiments in Plato and Augustine. And upon that note --- number and geometry as a path to the eternal ---  we ended our seminar.

\section{Mathematics, humanism and freedom} 
\begin{quote}
    {\em Wir m\"{u}ssen wissen – wir werden wissen.\/} ---Hilbert~\cite{hilbert_naturkennen}
\end{quote}

Mathematics, now, is in the midst of a perfect storm. The pressures on young mathematicians to conform, to take safe paths, have been increasing for many years, reflecting a hollowing out of our broader enterprise. But as large language models are carrying out mathematical research, and the social contract between societies and universities is fraying, mathematicians are called upon to account for the human value of our profession. We must take stock of the fact that while we define ourselves through our creativity, our discipline has been defined by others for its utility. Broader public narratives leave us trapped in such tropes as lonely genius, madness and mysticism.

Yet a look at the longer past provides hope. Mathematics is resilient---both the ideas and the people that create them. The modern profession of mathematics reflects our time, place and society. But it swims in a broader stream of thought, which, as we have seen in  the previous section, transcends disciplines and nations. This broader stream, which we sought to explore in our seminar, offers different tools to question the fundamental nature of mathematical thought and many visions for our collective future.

Mathematics is the product of free minds responding to an ever-changing world, and its history demonstrates, again and again, an ability to escape our best laid plans. Our very first reading of the seminar was Hilbert’s 1930 radio address in K{\"o}nigsberg, concluding with his stirring words above. This address came a day after a young Austrian named Kurt G{\"o}del gave a lecture that would reshape the foundations of the subject, and nine years before a war that would transform the world. What endures in mathematics lies not in its institutions or concepts or even theorems, but in its relationship to the human mind. It is at this juncture that the true riches of our subject reside; it is here that our seminar flourished; and   here where our discipline can find inspiration for its future, as it has so often in the past.

\appendix

\section{Reading List}
The ``large'' themes
covered in the seminar, together with the primary readings, were as follows.  The themes
were not covered consecutively -- we would often interweave them from one week to the next, which, we felt, created a somewhat lighter spirit;
and although we largely present the readings in chronological order below, they were usually not covered in that way.

\begin{itemize}

\item {\em The emergence of  modern geometry over the 19th century (5 sessions)\/}  

We read  Lobachevsky's work~\cite{lobachevsky_parallells}  on hyperbolic geometry,   the work of Riemann~\cite{riemann_hypothesen} and
Helmholtz~\cite{helmholtz_thatsachen,helmholtz_origin,helmholtz_origin_part2},    and concluded with Klein ~\cite{klein_nicht_euklidische,klein_erlangen}. We also read Hilbert's axiomatization~\cite{hilbert_grundlagen} of Euclidean geometry together with Brouwer's  inaugural address ~\cite{brouwer_intuitionism}.
We concluding by reading Mumford~\cite{mumford_stochasticity} and  Thurston~\cite{thurston_proof_progress,thurston_weeks_three_dimensional}
for some modern viewpoints. Non-technical material included
 correspondence between Gauss and Bolyai Sr;  Poincar{\'e}'s review of Hilbert's book;
and popular writings of Einstein.

\item {\em From atomism to probability (2 session) \/} 

We traced a path  back from Brownian motions to arguments about the molecular nature of matter, reading Maxwell~\cite{maxwell_molecules,maxwell_gases} and     
 excerpts from Mach~\cite{mach_mechanics},    as well as
Einstein~\cite{einstein_brownian},  and parts of Wiener~\cite{wiener_differential_space}.
Other material included 
 Wiener’s writings on science
  and society and Perrin’s experimental verification of Einstein’s work.

\item {\em From laws of thought to the automation of mathematics (3 sessions) \/}  

We studied  ``classical'' theories relating thought and logic, including
Boole~\cite{boole_laws}, Leibniz~\cite{leibniz_characteristic}, Peirce~\cite{peirce_algebra_writings}, 
and material from the dawn of the computer age --  McCulloch and Pitts ~\cite{mcculloch_pitts_calculus}
and von Neumann~\cite{vonneumann_mathematician,vonneumann_automata,vonneumann_computer_brain}.
We then read Dick's sociological study~\cite{dick_aftermath}  of early automated theorem proving   
and contrasted it with recent examples chosen by Constantin Kogler:   Georgiev {\em et al}~\cite{georgiev_mathematical}, Bryan {\em et al}~\cite{bryan_motivic}, Schmitt~\cite{schmitt_extremal}, Sothanaphan~\cite{sothanaphan_erdos}.  
Other material included Turing \cite{Turing}, more Leibniz, and Hilbert's K{\"o}nigsberg address.

\item {\em Newton and 17th century science (4 sessions) \/}.

  We studied the first three chapters of the {\em Principia\/}~\cite{newton_principia} carefully,
  culminating in the derivation of the universal law of gravitation from Kepler's area law. 
We studied also Newton's treatment of fluxions~\cite{newton_fluxions},  which we contrasted with Leibniz \cite{leibniz_new_method},
and  Newton's early adoption, and later rejection, of Descartes' work~\cite{descartes_geometrie}.
Finally we went back to Kepler~\cite{kepler_harmonice}. 
Other material included works both about and by Newton \cite{guicciardini, keynes_newton, newton_prophecies}, modern treatments
of the Kepler problem, and
Pauli's~\cite{pauli_kepler} Jung-ian analysis of Kepler.
\end{itemize}

We also had a number of standalone sessions, not part of any planned arcs of the seminar -- although in some cases
they reflected ideas that occurred across many of the other lectures.

\begin{itemize}
\item {\em Recreational mathematics. \/}   We studied recreational mathematics through
readings curated with Evelyn Lamb: Martin Gardner's columns on tiling problems~\cite{gardner_tessellating},
  Schattschneider~\cite{Schattschneider,schattschneider_tiling} and  Wang~\cite{wang_games}.

\item {\em The birth of category theory.\/} 
We studied the founding papers of category theory -- work of Eilenberg, Maclane and Steenrod~\cite{eilenberg_maclane_homomorphisms,eilenberg_maclane_extensions,eilenberg_maclane_equivalences, eilenberg_steenrod_axiomatic}, alongside the collected correspondence between Steenrod and Raymond Wilder~\cite{wilder_steenrod_letters}. These readings were chosen with Alma Steingart.

\item {\em Women in mathematics.\/}
 The creation of opportunities for women in mathematics was discussed through the lens of  broader social movements.
 The readings (curated with Karen Uhlenbeck) were primarily biographical:  
 Kessel~\cite{kessel_tenured_women}, Reid~\cite{reid_julia_robinson}, Jackson~\cite{jackson_miracle_group}, Uhlenbeck~\cite{uhlenbeck_glass_ceiling}.

\item {\em Mathematics and divinity.\/}  We looked at the relation of  mathematics and theology across several  eras: from classical Greece,  
Plato~\cite{plato_epinomis}; from early Christianity, Augustine of Hippo~\cite{augustine_free_choice}, and from the 20th century,
the wartime correspondence between Andr\'{e} and Simone Weil
~\cite{weil_correspondence}. This reading was curated with David Nirenberg, whose book~\cite{nirenberg_uncountable} was also part of the reading.

 \end{itemize}

\section{Statement on AI usage}
The authors made use of large language models  to critique and analyze the text, as well as to assist with preparation and verification of   some items in the bibliography.

\bibliographystyle{amsplain}
\bibliography{folklore}

\end{document}